\documentclass{elsarticle}

\def\opn#1#2{\def#1{\operatorname{#2}}} 
\opn\con{conv}
\opn\gr{gr}

\let\to=\rightarrow

\def\Implies{\ifmmode\Longrightarrow \else
     \unskip${}\Longrightarrow{}$\ignorespaces\fi}
\def\implies{\ifmmode\Rightarrow \else
     \unskip${}\Rightarrow{}$\ignorespaces\fi}
\def\iff{\ifmmode\Longleftrightarrow \else
     \unskip${}\Longleftrightarrow{}$\ignorespaces\fi}

\let\:=\colon
\let\epsilon=\varepsilon
\let\phi=\varphi

\usepackage[english]{babel}
\usepackage{amssymb}
\usepackage{amsthm}
\usepackage{amsfonts}

\def\dj{d\kern-0.4em\char"16\kern-0.1em}
\def\Dj{\mbox{\raise0.3ex\hbox{-}\kern-0.4em D}}
\newfont{\cir}{wncyr10 scaled 1095}
\newfont{\ciri}{wncyi10 scaled 1095}
\newfont{\cirsc}{wncysc10 scaled 1095}

\newtheorem{definition}{Definition}[section]
\newtheorem{theorem}{Theorem}[section]

\newtheorem{example}{Example}[section]

\newtheorem{corollary}{Corollary}[section]

\newcommand{\eqnsection}{
\renewcommand{\theequation}{{\thesection.\arabic{equation}}}
\renewcommand{\con}{\mathrm{conv}}
\makeatletter \csname @addtoreset\endcsname{equation}{section}
\makeatother} \eqnsection
\def\Proof{\removelastskip\vskip\baselineskip\noindent{\it \textbf{Proof}.\quad}\ignorespaces}
\title{Some new common fixed point results of generalized contractive
non-self multi-valued mappings}

\author{Farshid Khojasteh $^{a,1}$ and Vladimir
Rako\v cevi\' c $^{b,2}$}

\fntext[fn1] {$^*$Corresponding Author}

\fntext[fn2] {Supported  by Grant No. 174025 of the Ministry of
Science, Technology and Development, Republic of Serbia.}

\address[rvt]{
Department of Mathematics, Islamic Azad University, Arak-Branch,
Arak, Iran, f-khojaste$@$iau-arak.ac.ir }

\address[rvt]{University of Ni\v s, Faculty of Sciences and
Mathematics, Department of Mathematics, Vi\v segradska 33, 18000
Ni\v s, Serbia, vrakoc@ptt.rs }

\begin{document}

\begin{frontmatter}
\begin{abstract}
In this paper, we study some new common fixed point results of
generalized contractive non-self multi-valued mappings on complete
metric space.

\end{abstract}


\end{frontmatter}

\section{Introduction and Preliminaries
}

Many applications of the contraction mapping theorems occur in a convex setting and the
mapping involved is not necessarily self-mapping.
In this paper  we introduce and study  some new common fixed point results of generalized
contractive non-self multi-valued mappings on complete metric space.
We obtain a necessary and sufficient
condition to exist a common endpoint for such type of mappings, and
we apply our main results to obtain some common fixed point results  for
multi-valued mappings and
for single-valued mappings.
Our results are related to the well known results of Nadler
\cite{r11} and \' Ciri\' c \cite{r14}, and to the recent results of
Amini-Harandi  \cite{r8} and Moradi and Khojasteh \cite{r16}.

Let $(X,d)$ be a metric space, and let  $P_{cl,bd}(X)$ denote the
classes of all non empty, closed and bounded subsets of $X$. Let
$T:X \to P_{cl,bd}(X)$ be a multi-valued mapping on $X$. A point
$x\in X$ is called a fixed point of $T$ if $x\in Tx$. Set
$Fix(T)=\{x\in X:x\in Tx\}$. An element $x\in X$ is said to be an
{\bf strict fixed point}(endpoint) of $T$, if $Tx=\{x\}$.

Note that {\bf ``strict fixed point"} is more natural
than {\bf ``endpoint"} and we apply ``strict fixed point"
instead of ``endpoint", throughout this paper
(see \cite{v3} and \cite{v4} for more details).\\
The set of all strict fixed points of $T$
denotes by $SFix(T)$. Obviously, $SFix(T)\subseteq Fix(T)$. The famous
theorem is doe to Nadler \cite{r11}. He extended the Banach
contraction principle to multi-valued mappings. Many authors have
studied the existence and uniqueness of strict fixed points for a multi-valued
mappings in metric spaces, see for example (\cite{r8},\cite{v1} and \cite{c1}-\cite{r11}) and references therein.


Let  $H$ be
 the Hausdorff metric on $P_{cl,bd}(X)$ induced by $d$,
that is,
\begin{equation}\label{d3}
H(A,B):=\max\Big{\{}\sup_{x\in B}d(x,A),\sup_{x\in A}d(x,B)\Big{\}},
\quad A,B\in P_{cl,bd}(X).
\end{equation}

A mapping $T:X\to P_{cl,bd}(X)$ has the approximation strict fixed point
property \cite{r8}, if
\begin{equation}\label{d9}
\inf_{x\in X}\sup_{y\in Tx} d(x,y)=0.
\end{equation}
If $T:X\to X$ is  a single-valued mapping, then $T$ has the
approximate strict fixed point property if and only if $T$ has the approximate
fixed point property, i.e.,
$
\inf_{x\in X}d(x,Tx)=0.
$

Recently, in  2010,  Amini-Harandi  \cite{r8} proved the following
very interesting fixed point theorem
\begin{theorem}\label{fdq3}
Let $(X,d)$ be a complete metric space. Suppose that $T:X\to
P_{cl,bd}(X)$ is a multi-valued mapping that satisfies
\begin{equation}\label{pp3}
H(Tx,Ty)\leq\psi(d(x,y)),
\end{equation}
for each $x,y\in X$, where $\psi:[0,+\infty)\to [0,+\infty)$ is
upper semicontinuous, $\psi(t)<t$ for all $t>0$, and satisfies
$\liminf_{t\to\infty}(t-\psi(t))>0$. Then $T$ has a unique strict fixed point
if and only if $T$ has the approximate strict fixed point property.
\end{theorem}

In 1971, \' Ciri\' c \cite{r14}, among other things, introduced and
proved the following interesting generalization of  of contraction
mapping principle.

\begin{theorem}
 Let $(X,d)$ be a complete  metric space. Let $f$ be a general
contractive mapping on $X,$ that is, there exists $\lambda \in
[0,1)$ satisfying
\[
d(fx,fy)\leq \lambda \cdot \max \bigg{\{}d(
x,y),d(fx,x),d(fy,y),\frac{%
d(fx,y)+d(fy,x)}2\bigg{\}},
\]
for all $x,y\in X.$  Then $f$ has a unique fixed point $u\in X$.
Moreover, for every $x\in X$, $ u=\lim _{n}f^nx. $
\end{theorem}

Very recently,  in 2011, Moradi and Khojasteh \cite{r16} proved the
next result as an extension of Theorem \ref{fdq3} when $\psi$ is a
nondecreasing map:
\begin{theorem}\label{fdq}
Let $(X,d)$ be a complete metric space. Suppose that $T:X\to
P_{cl,bd}(X)$ is a multi-valued mapping that satisfies
\begin{equation}\label{pp}
H(Tx,Ty) \leq\psi
\bigg(\max\Big{\{}d(x,y),d(x,Tx),d(y,Ty),\frac{d(x,Ty)+d(y,Tx)}{2}\Big{\}}
\bigg),
\end{equation}
for each $x,y\in X$, where $\psi:[0,+\infty)\to [0,+\infty)$ is
upper semicontinuous, $\psi(t)<t$ for all $t>0$, and satisfies
$\liminf_{t\to\infty}(t-\psi(t))>0$. Then $T$ has a unique strict fixed point
if and only if $T$ has the approximate strict fixed point property.
\end{theorem}


\section{Main results}

In this section we discuss our main result. Let $K$  be a closed
subset of $X$ and $\partial K$
 the boundary of $K$.
\begin{definition}
A non-self mapping $T:K\to P_{cl,bd}(X)$ has the approximate
$K-$boundary strict fixed point property, if
\begin{equation}\label{d9}
\inf\{H(\{x\},Tx):x\in\partial K\}=0.
\end{equation}
Also, we say that two non-self multi-valued mappings $T,S:K\to
P_{cl,bd}(X)$ have common approximate $K-$boundary strict fixed point property
if there exist a sequence $\{x_n\}\subset\partial K$ such that

\begin{equation}\label{d229}
\lim_{n\to\infty}H(\{x_n\},Tx_n)=0 \ \ \ and \ \ \
\lim_{n\to\infty}H(\{x_n\},Sx_n)=0.
\end{equation}
Note that if $T,S$ be two single-valued mappings from $K$ into $X$,
then have common approximate $K-$boundary strict fixed point property if and
only if $T$ and $S$ have common approximate fixed point property,
i.e. there exists a sequence $\{x_n\}\subset\partial K$ such that
\begin{equation}\label{dg229}
\lim_{n\to\infty}d(x_n,Tx_n)=0 \ \ \ and \ \ \
\lim_{n\to\infty}d(x_n,Sx_n)=0.
\end{equation}
\end{definition}
\begin{theorem}\label{ww2}
Let $(X,d)$ be a complete metric space and $K$ be a closed subset of
$X$. Suppose that $T,S:K\to P_{cl,bd}(X)$ are two multi-valued
mappings such that
\begin{equation}\label{pp}
H(Tx,Sy)\leq\psi \biggl (N_{T,S}(x,y)\biggr ),
\end{equation}
for each $x,y\in X$, where
\begin{equation}
N_{T,S}(x,y)=\max\biggl{\{}d(x,y),d(x,Tx),d(y,Sy),\frac{d(y,Tx)+d(x,Sy)}{2}\biggr{\}},
\end{equation}
and $\psi:[0,+\infty)\to [0,+\infty)$ is upper semi-continuous,
$\psi(t)<t$ for all $t>0$, and satisfies
$\liminf_{t\to\infty}(t-\psi(t))>0$. Then $T,S$ have a unique common
strict fixed point in $K$ if and only if $T,S$ have the common approximate
$K-$boundary strict fixed point property. Also, $SFix(T)=Fix(T)=Fix(S)=SFix(S)$.
\end{theorem}
\Proof It is clear that if $T$ and $S$ have common strict fixed point, then
$T$ and $S$ has the approximate $K-$boundary strict fixed point property.
Conversely, suppose that $T$ and $S$ have the common approximate
$K-$boundary strict fixed point property; then there exists a sequence
$\{x_n\}\subset\partial K$ such that ${\lim_n}H(\{x_n\},Tx_n)=0$ and
${\lim_n}H(\{x_n\},Sx_n)=0$. For all $m,n\in\Bbb  N$ we have
\begin{equation}\label{fffg}
\begin{array}{lll}
N_{T,S}(x_n,x_m) &\leq& \max \biggl
\{d(x_n,x_m),H(\{x_n\},Tx_n),H(\{x_m\},Sx_m),
\frac{H(\{x_n\},Tx_m)+H(\{x_m\},Sx_n)}{2}\biggr \}\\
&\leq& d(x_n,x_m)+H(\{x_n\},Tx_n)+H(\{x_m\},Sx_m)\\
&=&d(x_n,x_m)-H(\{x_n\},Tx_n)-H(\{x_m\},Sx_m)\\
&+&2H(\{x_n\},Tx_n)+2H(\{x_m\},Sx_m)\\
&\leq&H(Tx_n,Sx_m)+2H(\{x_n\},Tx_n)+2H(\{x_m\},Sx_m)\\
&\leq&\psi(N(x_n,x_m))+2H(\{x_n\},Tx_n)+2H(\{x_m\},Sx_m).
\end{array}
\end{equation}
Since $\psi$ is u.s.c., $\psi(t)<t$ for all $t>0$ and
$\liminf_{t\to\infty}(t-\psi(t))>0$, from (\ref{fffg}), we have
\[
\limsup_{m,n\to\infty}N_{T,S}(x_n,x_m)=0.
\]
Thus, $\{x_n\}$ is a Cauchy sequence in $\partial K$. So there
exists $x_0 \in\partial K$ such that ${\lim_n}x_n=x_0$. Note that
if $N_{T,S}(x_{n_0},x_0)=0$ for some $n_0\in\Bbb N$, then
$x_0=x_{n_0}$, $d(x_0,Sx_0)=0$ and $d(x_0,Tx_0)=0$. This means
that $x_0\in \overline{Tx_0} = Tx_0\subset K$ and $x_0\in
\overline{Sx_0} = Sx_0\subset K$ this completes the proof. Thus,
suppose that $N_{T,S}(x_n,x_0)\Bbb \neq 0$ for all $n\in \Bbb N$.
We claim that $x_0\in Tx_0$ and $x_0\in Sx_0$.\\
If $x_0 \Bbb  Notin Sx_0$ then $d(x_0,Sx_0) > 0$. For all $n \in
\Bbb N$,
\begin{eqnarray}\label{dg44}
\frac{d(x_n,Sx_0)+d(x_0,Tx_n)}{2} \leq
\frac{d(x_n,x_0)+d(x_0,Sx_0)+d(x_n,x_0)+d(x_n,Tx_n)}{2}.
\end{eqnarray}
Since ${\lim_n}d(x_{n},x_0)={\lim_n}d(x_{n},Tx_n)=0$ , $d(x_0,Sx_0)
> 0$ and (\ref{dg44}) holds, there exists $N_1 \in \Bbb N$ such that
for all $n \geq N_1$
\begin{eqnarray}\label{d46}
\frac{d(x_n,Sx_0)+d(x_0,Tx_n)}{2} < d(x_0,Sx_0).
\end{eqnarray}
So for all $n \geq N_1$, $N_{T,S}(x_n,x_0)=d(x_0,Sx_0)$ and hence
for all $n \geq N_1$, we have
\begin{equation}\label{fffg2}
\begin{array}{lll}
H(\{x_n\},Sx_0)-H(\{x_n\},Tx_n)&\leq&H(Tx_n,Tx_0)\\
&\leq&\psi(N_{T,S}(x_n,x_0))\\
&<&N_{T,S}(x_n,x_0)\\
&=&d(x_0,Sx_0).
\end{array}
\end{equation}
Thus for each $n\geq N_1$, we have
\begin{equation}
H(\{x_n\},Sx_0)-H(\{x_n\},Tx_n)<d(x_0,Sx_0).
\end{equation}
Therefore, $H(\{x_0\},Sx_0)\leq d(x_0,Sx_0)$ and this is a
contradiction. Hence, $x_0\in Sx_0$. The same argument implies
$x_0\in Tx_0$.
Now we prove that $Sx_0=\{x_0\}$(one can apply similar argument for $T$).\\
Suppose that $H(\{x_0\},Sx_0)\Bbb  \neq 0$. For all $n \in \Bbb N$,
\begin{equation}\label{d4f8}
\begin{array}{lll}
H(\{x_n\},Sx_0)-H(\{x_{n}\},Tx_{n}) &\leq& H(Tx_{n},Sx_0)\\
&\leq& \psi(N_{T,S}(x_{n},x_0))\\
&\leq &N_{T,S}(x_{n},x_0)\\
&\leq& d(x_n,x_0)+H(\{x_0\},Sx_0)+H(\{x_n\},Tx_n),
\end{array}
\end{equation}
We conclude from (\ref{d4f8}) that
${\lim_n}N_{T,S}(x_n,x_0)=H(\{x_0\},Sx_0)$. Since $\psi$ is l.s.c
\begin{equation}
\limsup_n\psi(N_{T,S}(x_n,x_0))\leq \psi(H(\{x_0\},Sx_0)).
\end{equation}
Also, we obtain from (\ref{d4f8}) that
${\lim_n}\psi(N_{T,S}(x_n,x_0))=H(\{x_0\},Sx_0)$. Therefore,
\begin{equation}
H(\{x_0\},Sx_0)\leq \psi(H(\{x_0\},Sx_0))
\end{equation}
and this is a contradiction. Hence $Sx_0=x_0$.\\
Now, we claim that such common strict fixed point is unique. In other words,
$SFix(T)=Fix(T)=Fix(S)=SFix(S)$.  Let $y\in Fix(T)$ be arbitrary. We
need to show that $y=x_0$. If $y\neq x_0$ then
\begin{equation}\label{fdg66}
d(x_0,y)\leq H(\{x_0\},Ty)=H(Sx_0,Ty)\leq
\psi(N_{T,S}(x_0,y))<N_{T,S}(x_0,y),
\end{equation}
where $N_{T,S}(x_0,y)=\max\{d(x_0,y),[{d(x_0,Ty)+d(y,x_0)}]/{2}\}$.
Since $y\in Ty$, $d(x_0,Ty)\leq d(x_0,y)$ and hence
$N_{T,S}(x_0,y)=d(x_0,y)$. From (\ref{fdg66}) we conclude that
$d(x_0,y)<N_{T,S}(x_0,y)$ and this is a contradiction. Therefore,
$SFix(T)=Fix(T)$.\quad$\Box$

\begin{example}
Let $X=\{x,y,z,w\}$, $K=\{x,y\}$ and let
\[
\begin{array}{lll}
d(x,y)=\frac{3}{2}&d(z,y)=\frac{6}{5}& d(w,y)=\frac{13}{10},\\
d(x,z)=1 & d(x,w)=1 & d(z,w)=1,\\
d(a,a)=0 & \forall \ a\in X,&\\
d(a,b)=d(b,a) & \forall \ a,b\in X &.
\end{array}
\]
Precisely, $(X,d)$ is a complete metric space and $\partial K=K$.
Now let $T,S:K\to P_{cl,bd}(X)$ are defined by
$Tx=\{y,z\}, Ty=\{y\}$ and $Sx=\{z,w\}, Sy=\{y\}$.\\
If we define
\begin{equation}
\psi(t)=\left\{
\begin{array}{lll}
\frac{14}{15}t&& 0\leq t\leq \frac{3}{2}\\\\
\frac{t}{2}t&& t>\frac{3}{2}
\end{array}
\right.
\end{equation}
then $\psi(t)<t$, $\liminf_{t\to\infty}(t-\psi(t))=+\infty$ and $\psi$ is lower semi continuous.
We have also,
\[
H(Tx,Sy)=d(z,y)=\frac{6}{5} \ \ \ and \ \ \ H(Ty,Sx)=\max\{d(y,z),d(y,w)\}=\frac{13}{10}
\]
One can show that $T$ satisfies in all conditions of Theorem \ref{ww2}.\\
On the other hand, S(y)=T(y)=\{y\} is only common strict fixed point of $T,S$ and if we consider the sequence $\{x_n\}_{n\in\mathbb{N}}=\{y\}$.
Then,
$\{x_n\}\subset \partial K$ and we have
\[
\lim_{n\to\infty}H(\{x_n\},Tx_n)=0 \ \ \ and \ \ \
\lim_{n\to\infty}H(\{x_n\},Sx_n)=0.
\]
\end{example}

\medskip The investigations of the contractive condition in  the next
theorem have been motivated by the results of \cite{c1}.

\begin{theorem}\label{kkl}
Let $(X,d)$ be a complete metric space and $K$ be a nonempty closed
subset of $X$. Let $T,S:K\to P_{cl,bd}(X)$ be two multi-valued
mappings such that
\begin{equation}\label{ffdgg}
\begin{array}{lll}
H(Tx,Sy)&\leq &\alpha d(x,y)+\beta\max\{d(x,Tx),d(y,Sy)\}\\
&&+\gamma\max\{d(x,Tx)+d(y,Sy),d(x,Ty)+d(y,Sx)\},
\end{array}
\end{equation}
for all $x,y\in K$, where $\alpha,\beta,\gamma\geq 0$ are such that
\begin{equation}
\alpha+2\gamma+\beta<1
\end{equation}
Then $T,S$ have a unique common strict fixed point in $K$ if and only if $T,S$
have the common approximate $K-$boundary strict fixed point property. Also,
$SFix(T)=Fix(T)=Fix(S)=SFix(S)$.
\end{theorem}
\Proof If is clear that if $T$ and $S$ have common strict fixed point, then
$T$ and $S$ has the approximate $K-$boundary strict fixed point property.
Conversely, suppose that $T$ and $S$ have the common approximate
$K-$boundary strict fixed point property; then there exists a sequence
$\{x_n\}\subset\partial K$ such that ${\lim_n}H(\{x_n\},Tx_n)$ $=0$
and ${\lim_n}H(\{x_n\},Sx_n)=0$. For all $m,n\in\Bbb N$ we have
\begin{equation}\label{gffg}
\begin{array}{lll}
d(x_n,x_m)&\leq& H(\{x_n\},Tx_n)+H(\{x_m\},Sx_m)+H(Tx_n,Sx_m)\\
&\leq&H(\{x_n\},Tx_n)+H(\{x_m\},Sx_m)+\alpha d(x_n,x_m)\\
&+&\beta\max\{d(x_n,Tx_n),d(x_m,Sx_m)\}\\
&+&\gamma\max\{d(x_n,Tx_n)+d(x_m,Sx_m),d(x_n,Tx_m)+d(x_m,Sx_n)\}\\
&\leq& H(\{x_n\},Tx_n)+H(\{x_m\},Sx_m)+\alpha d(x_n,x_m)\\
&+&\beta\max\{H(\{x_n\},Tx_n),H(\{x_m\},Sx_m)\}\\
&+&\gamma\max\{H(\{x_n\},Tx_n)+H(\{x_m\},Sx_m),H(\{x_n\},Tx_m)+H(\{x_m\},Sx_n)\}\\
&\leq&H(\{x_n\},Tx_n)+H(\{x_m\},Sx_m)+\alpha d(x_n,x_m)\\
&+&\beta(H(\{x_n\},Tx_n)+H(\{x_m\},Sx_m))\\
&+&\gamma(H(\{x_n\},Tx_n)+H(\{x_m\},Sx_m)+H(\{x_n\},Tx_m)+H(\{x_m\},Sx_n))\\
&\leq&H(\{x_n\},Tx_n)+H(\{x_m\},Sx_m)+\alpha d(x_n,x_m)\\
&+&\beta(H(\{x_n\},Tx_n)+H(\{x_m\},Sx_m))\\
&+&\gamma(2d(x_n,x_m)+H(\{x_n\},Tx_n)+H(\{x_m\},Sx_m)\\
&+&H(\{x_n\},Tx_n)+H(\{x_m\},Sx_m))\\
&=&(\beta+3\gamma+1)[H(\{x_n\},Tx_n)+H(\{x_m\},Sx_m))]+(\alpha+2\gamma)d(x_n,x_m).
\end{array}
\end{equation}
It means that,
\begin{equation}
d(x_n,x_m)\leq
\frac{\beta+3\gamma+1}{1-\alpha-2\gamma}[H(\{x_n\},Tx_n)+H(\{x_m\},Sx_m))]
\end{equation}
Thus ${\limsup}_{n,m\to\infty}d(x_n,x_m)=0$, i.e., $\{x_n\}$ is a
Cauchy sequence in $\partial K$ and so convergent to $x_0\in K$. We
claim that $Tx_0=\{x_0\}$ and $Sx_0=\{x_0\}$. If $Tx_0\neq\{x_0\}$
then $H(\{x_0\},Tx_0)\neq 0$. For all $n \in \Bbb  N$,
\begin{equation}\label{dg44t}
d(x_n,Tx_0)+d(x_0,Sx_n) \leq
2d(x_n,x_0)+H(\{x_0\},Tx_0)+H(\{x_n\},Sx_n).
\end{equation}
Since ${\lim}_n d(x_{n},x_0)={\lim}_n H(\{x_n\},Tx_n)=0$,
$d(x_0,Tx_0)
> 0$ and (\ref{dg44t}) holds, there exists $N_1 \in \Bbb  N$ such
that for all $n \geq N_1$
\begin{eqnarray}\label{d46h}
d(x_n,Tx_0)+d(x_0,Tx_n) < H(\{x_0\},Tx_0).
\end{eqnarray}
(one can apply similar argument for $S$). We have\\
\begin{equation}\label{d4f8h}
\begin{array}{lll}
H(\{x_n\},Tx_0)-H(\{x_{n}\},Sx_{n}) &\leq& H(Sx_{n},Tx_0)\\
&\leq& \alpha d(x_0,x_n)+\beta\max\{d(x_0,Tx_0),d(x_n,Sx_n)\}\\
& + & \gamma\max\{d(x_0,Tx_0)+d(x_n,Sx_n),d(x_n,Tx_0)+d(x_0,Sx_n)\\
&\leq&\alpha d(x_0,x_n)+\beta\max\{H(\{x_0\},Tx_0),H(\{x_n\},Sx_n)\}\\
& + & \gamma\max\{H(\{x_0\},Tx_0)+H(\{x_n\},Sx_n),H(\{x_n\},Tx_0)\\
&+&H(\{x_0\},Sx_n)
\end{array}
\end{equation}
Hence,  by taking limit from two sides of (\ref{d4f8h}) we conclude
that
\begin{equation}
H(\{x_0\},Tx_0)\leq(\beta+\gamma)H(\{x_0\},Tx_0).
\end{equation}
Since $\beta+\gamma<1$ we have a contradiction. Thus $Tx_0=\{x_0\}$.\\
Also, we claim that such common strict fixed point is unique. In other words,
$SFix(T)=Fix(T)=Fix(S)=SFix(S)$.  Let $y\in Fix(T)$ be arbitrary. We
need to show that $y=x_0$. If $y\neq x_0$ then
\begin{equation}\label{fdg66h}
\begin{array}{lll}
d(x_0,y)&\leq &H(\{x_0\},Ty)=H(Sx_0,Ty)\\
&\leq&\alpha d(x_0,y)+\beta\max\{d(y,Ty),d(x_0,Sx_0)\}\\
&+&\gamma\max\{d(y,Ty)+d(x_0,Sx_0),d(x_0,Ty)+d(y,Sx_0)\\
&=&(\alpha+2\gamma)d(x_0,y)\\
&<& d(x_0,y),
\end{array}
\end{equation}
and this is a contradiction. Therefore, $SFix(T)=Fix(T)$.\quad$\Box$
\begin{example}
Let $X=[0,2]$ endowed with the Euclidean metric, $K=\{0,1,2\}$ and $T,S:K\to P_{cl,bd}(X)$ are defined by
\[
T0=\{1,\frac{2}{5}\},~~T1=\{1\},~~ T2=\{\frac{2}{3}\}
\]
and
\[
S0=\{\frac{3}{4}\},~~ S1=\{1\},~~ S2=\{\frac{2}{5}\}.
\]
It is easily seen that
\[
\begin{array}{lll}
H(T0,S1)=\frac{3}{5} &H(T1,S0)=\frac{1}{4}&H(T0,S2)=\frac{3}{5}\\\\
H(T2,S0)=\frac{1}{12}&H(T1,S2)=\frac{3}{5}&H(T2,S1)=\frac{1}{3}.
\end{array}
\]
If we consider $\alpha=\beta=\frac{1}{3}$ and $\gamma=\frac{3}{25}$, then $\alpha+\beta+2\gamma=\frac{68}{75}<1$ and for
each $x,y\in K$ we have
\[
\begin{array}{lll}
H(Tx,Sy)&\leq& \frac{1}{3}d(x,y)+\frac{1}{3}\max\{d(x,Tx),d(y,Sy)\}\\\\
&+&\frac{3}{25}\max\{d(x,Tx)+d(y,Sy),d(x,Ty)+d(y,Sx)\}
\end{array}
\]
For example if $x=0$ and $y=1$ then $d(0,T0)=\frac{2}{5}$, $d(1,S1)=0$, $d(0,T1)=1$, $d(1,S0)=\frac{1}{4}$ and thus
\[
H(T0,S1)=\frac{3}{5}<\frac{37}{60}
=\frac{1}{3}+\frac{1}{3}\max\{\frac{2}{5},0\}+\frac{3}{25}\max\{\frac{2}{5},\frac{5}{4}\}.
\]
Therefore, $T$ satisfies in all conditions of Theorem \ref{kkl}.\\
On the other hand, S(1)=T(1)=\{1\} is only common strict fixed point of $T,S$ and if we consider the sequence $\{x_n\}_{n\in\mathbb{N}}=\{1\}$.
Then,
$\{x_n\}\subset \partial K$ and we have
\[
\lim_{n\to\infty}H(\{x_n\},Tx_n)=0 \ \ \ and \ \ \
\lim_{n\to\infty}H(\{x_n\},Sx_n)=0.
\]
\end{example}
\begin{corollary}\label{www2}
Let $(X,d)$ be a complete metric space and $K$ be a closed subset of
$X$. Suppose that $f,g:K\to X$ are two single-valued mappings such
that
\begin{equation}\label{epp}
d(fx,gy)\leq\psi \biggl (N_{f,g}(x,y)\biggr ),
\end{equation}
for each $x,y\in X$, where $\psi:[0,+\infty)\to [0,+\infty)$ is
upper semi-continuous, $\psi(t)<t$ for all $t>0$, and satisfies
$\liminf_{t\to\infty}(t-\psi(t))>0$. Then $f,g$ have a unique common
fixed point in $K$ if and only if $f,g$ have the common approximate
fixed point on $\partial K$.
\end{corollary}
\begin{proof}
If we take $Tx=\{f(x)\}$ and $Sx=\{g(x)\}$ and using Theorem
\ref{ww2}, we obtain desired result.\quad$\Box$
\end{proof}
\begin{corollary}\label{kww2}
Let $(X,d)$ be a complete metric space and $K$ be a closed subset of
$X$. Suppose that $f,g:K\to X$  are two single-valued mappings such
that $gf(\partial K)\subset \partial K$, $f(\partial K)\subset K$
and $f,g$ satisfy the following:
\begin{description}
\item[(A)] For each $x,y\in K$
\begin{equation}\label{eppk}
d(fx,gy)\leq\psi \biggl (N_{f,g}(x,y)\biggr ),
\end{equation}
\item[(B)] $d(x,gx)\leq d(x,fx)$.
\end{description}
where $\psi:[0,+\infty)\to [0,+\infty)$ is upper semi-continuous,
$\psi(t)<t$ for all $t>0$, and satisfies
$\liminf_{t\to\infty}(t-\psi(t))>0$. Then $f,g$ have the common
approximate fixed point on $\partial K$.
\end{corollary}
\begin{proof}
Let $x_0\in\partial K$ and let $x_1=fx_0\in K$. Also,
$x_2=gx_1=gfx_0\in\partial K$. Thus we can choose a sequence
$\{x_n\}\subset K$ such that $\{x_{2n}\}\subset\partial K$. We
have $x_{2n+2}=gx_{2n+1}$ and $x_{2n+1}=fx_{2n}$, for all
$n\in\Bbb  N$. We conclude that
\begin{equation}\label{h32}
\begin{array}{lll}
d(x_{2n+1},x_{2n+2})&=&d(fx_{2n},gx_{2n+1})\\
&\leq& \psi(d(x_{2n},x_{2n+1}))\\
&<&d(x_{2n},x_{2n+1}).
\end{array}
\end{equation}
Also,
\begin{equation}\label{h33}
\begin{array}{lll}
d(x_{2n},x_{2n+1})&=&d(fx_{2n},gx_{2n-1})\\
&\leq& \psi(d(x_{2n},x_{2n-1}))\\
&<&d(x_{2n},x_{2n-1}).
\end{array}
\end{equation}
Therefore, for each $k\in\Bbb  N$ we have $
d(x_{k},x_{k+1})<d(x_{k},x_{k-1}). $ Hence $\{d(x_n,x_{n+1})\}$
is a nondecreasing sequence an so is convergent to $r\geq 0$.
Assume that $r>0$; then by (\ref{h32}),(\ref{h33}) we obtain $
r\leq \psi(r), $ and this is a contradiction. It means that
$r=0$. In other words, $\lim_{n\to\infty}d(x_{2n},fx_{2n})=0$ and
by $(B)$ We have $\lim_{n\to\infty}d(x_{2n},gx_{2n})=0$. Thus,
$f$ and $g$ have common approximate fixed point property on
$\partial K$.
\end{proof}
\begin{corollary}\label{jww2}
Let $(X,d)$ be a complete metric space and $K$ be a closed subset of
$X$. Suppose that $f,g:K\to X$ are two single-valued mappings such
that $gf(\partial K)\subset \partial K$, $f(\partial K)\subset K$
and $f,g$ satisfy the following:
\begin{description}
\item[(A)] For each $x,y\in K$
\begin{equation}\label{jppk}
d(fx,gy)\leq\psi \biggl (N_{f,g}(x,y)\biggr ),
\end{equation}
\item[(B)] $d(x,gx)\leq d(x,fx)$.
\end{description}
where $\psi:[0,+\infty)\to [0,+\infty)$ is upper semi-continuous,
$\psi(t)<t$ for all $t>0$, and satisfies
$\liminf_{t\to\infty}(t-\psi(t))>0$. Then $f,g$ have the common
fixed point on $K$.
\end{corollary}
\begin{proof}
 Using Corollary \ref{kww2} and Corollary \ref{www2} we obtain
desired result.
\end{proof}
\section*{Acknowledgment}
We would like to thank the anonymous referees for their useful comments and suggestions.

\end{document}